\documentclass[12pt]{amsart}
\title[Asdim of one relator groups]{Asymptotic dimension of one relator groups}
\usepackage{amsmath,amssymb,amsthm,hyperref}

\DeclareMathOperator{\asdim}{asdim}

\DeclareMathOperator{\dist}{dist}
\theoremstyle{plain}
\newtheorem{theorem}{Theorem}[section]

\theoremstyle{remark}

\theoremstyle{definition}
\newtheorem{definition}[theorem]{Definition}
\theoremstyle{remark}

\author{Dmitry Matsnev}
\address{Department of Mathematics\\The Pennsylvania State University\\ University Park, PA 16802\\U.S.A.}
\curraddr{Departamento de Matem\'atica\\Instituto Superior T\'ecnico\\Av. Rovisco Pais\\1049-001 LISBOA\\Portugal}
\email{matsnev@math.ist.utl.pt}
\subjclass[2000]{Primary 20F65; Secondary 20F05, 20F69}
\keywords{asymptotic dimension, one relator groups}
\date{July 20, 2006}
\begin{document}
\begin{abstract}
 We show that one relator groups viewed as metric spaces with respect to the word-length metric have finite asymptotic dimension in the sense of Gromov and give an estimate of their asymptotic dimension in terms of the relator length.
\end{abstract}
\maketitle

\section{Introduction}

A finitely generated group has asymptotic dimension not more than $n$ if its underlying metric space (with respect to the length function corresponding to a given set of generators) has the following property: for any $R>0$ there exists a uniformly bounded cover of the group such that any $R$-ball meets not more than $(n+1)$ elements of the cover (see~\cite{Gromov}). There is a way to extend this notion to all countable groups, not necessarily finitely generated.

Not all groups have finite asymptotic dimension ($\mathbb Z \wr \mathbb Z$ is the standard example of a group which does not satisfy the definition above for any $n$), but for some classes of groups it is known that their asymptotic dimension is finite: for instance, Gromov noticed this for the case of hyperbolic groups in~\cite{Gromov} (for precise proof, see~\cite{RoeCoarseGeometry}), Ji proved this for arithmetic groups in~\cite{Ji}, and of course there are more examples.

The question whether a given group has finite asymptotic dimension has drawn more attention after results of Yu, who showed in~\cite{Yu} that such group satisfies the Novikov conjecture, and  Higson and Roe, who deduced the exactness of such group in~\cite{HigsonRoe}.

The second part of the title concerns one relator groups. These are the ones which admit a presentation $< S | r >$, where the generating set $S$ is at most countable, and the only relator $r$ is a word in $S$. Historically, one relator groups appeared as the fundamental groups of $2$-manifolds, and one may regard them as a simplest class of groups which are  close to free groups.

Within the scope of this work, the main motivation for the study of the asymptotic dimension of one relator groups were the results of Guentner (see~\cite{Guentner}), showing that one relator groups are exact, and the ones of Beguin, Bettaieb, and Valette in~\cite{BeguinBettaiebValette}, showing that the Baum-Connes conjecture (and therefore the Novikov conjecture as well) holds for such groups. A very natural question which arises immediately after comparing all the results mentioned above is: do one relator groups have finite asymptotic dimension? In this paper we are giving an affirmative answer to this question and estimate asymptotic dimension in terms of the relator length.

\section{Asymptotic dimension of groups}

First we outline a few basic definitions and some theorems on asymptotic dimension which shall be used later. For more detailed discussion and precise proofs, consult~\cite{RoeCoarseGeometry}.

\begin{definition}
Given a metric space $G$, we say that its \emph{asymptotic dimension} does not exceed $n$ and write $\asdim G\leq n$ if for any $R>0$ there exists a uniformly bounded cover of $G$ such that any $R$-ball in $G$ meets not more than $(n+1)$ elements of the cover. The \emph{asymptotic dimension} of $G$, $\asdim G$, is then the minimal $n$ satisfying this condition.
\end{definition}

Any finitely generated group $\Gamma$ with a generating set $S$ can be endowed with a word-length metric
\begin{equation*}
\dist(\gamma_1,  \gamma_2)=\mbox{length of a shortest word in $S\cup S^{-1}$ representing $\gamma_1^{-1}\gamma_2$}.
\end{equation*}
Since any two metrics for the same group, arising form different generating sets, are bi-Lipschitz equivalent, the notion of asymptotic dimension of a finitely generated group, viewed as a discrete metric space as stated above, is independent of the particular generating set.

\begin{definition}
For a (not necessarily finitely generated) group $G$, its asymptotic dimension, $\asdim G$, is defined to be a supremum of $\asdim\Gamma$ over all finitely generated subgroups $\Gamma$ of $G$.
\end{definition}
This definition is consistent, for any inclusion of a subgroup into the over-group is a coarse embedding.

Now we collect some facts on the asymptotic dimension of groups, to be used in our discussion later on.

\begin{theorem}\label{thmSubgroup}
For a subgroup $H$ of a group $G$, $\asdim H\leq\asdim G$.
\end{theorem}
\begin{theorem}[Bell and Dranishnikov,~\cite{BellDranishnikov}]\label{thmHNNextension}
For an HNN extension $*_AG$ of a group $G$, $\asdim *_AG\leq1+\asdim G$
\end{theorem}
The original theorem in~\cite{BellDranishnikov} was formulated for a finitely generated base group $G$, but since any finitely generated subgroup of the HNN extension of $G$ is a subgroup of an HNN extension of some finitely generated subgroup of $G$, one can run the original argument of Bell and Dranishnikov for every finitely generated subgroup of $*_AG$ to obtain the theorem in the form we formulated here.
\begin{theorem}[Bell and Dranishnikov,~\cite{BellDranishnikov}]\label{thmFreeProduct}
For a free product $G*H$ of two groups $G$ and $H$, $\asdim G*H\leq\max\{\asdim G, \asdim H, 1\}$.
\end{theorem}
Again, the assumption that both $G$ and $H$ are finitely generated is not crucial for the proof.

Finally, as a base for our inductive arguments in the next section, we state that the asymptotic dimension of a finite group is $0$, and the asymptotic dimension of a free group is $1$.

\section{One relator groups}

Throughout this section, let $G$ be a one relator group with a (possibly infinite) generating set $S$ and relator $r$, that is, $G=<S|r>$ is a quotient of a free group on $S$ by the minimal normal subgroup generated by $r$. We assume that $r$, a finite word in $S\cup S^{-1}$, is cyclically reduced as a word in the free group on $S$, and use $|r|$ to denote its length in this free group. To omit the trivial cases, we assume that $S$ contains at least two elements and $|r|>0$.

For any real number $x$, we denote by $\lceil x\rceil$ the minimal integer greater or equal to $x$.

\begin{theorem}
In the notations above, $\asdim G\leq\lceil\frac{|r|}2\rceil$.
\end{theorem}

The rest of this section is devoted to the proof of this theorem.

First note that we can assume that $G$ is finitely generated and every letter of $S$ appears in $r$. Indeed, $G$ is isomorphic to a free product of a finitely generated one relator group $\Gamma$ with relator $r$ and generating set consisting of letters which appear in $r$ and the free group on all other letters. According to Theorem~\ref{thmFreeProduct}, $\asdim G\leq\max\{\asdim\Gamma, 1\}$. If we can prove that $\asdim\Gamma\leq\lceil\frac{|r|}2\rceil$, the statement of the  theorem for $G$ will follow.

The argument is based on the induction on the length of $r$. For $|r|=1$ the group $G$ is isomorphic to a free group on all letters in $S$ except the one which appears in $r$. Thus $\asdim G=1\leq1=\lceil\frac{|r|}2\rceil$. 

For the inductive step suppose that the statement of the theorem has been proven for all one relator groups with relator length strictly less than $|r|$. Following the standard arguments of Magnus and Molchanovskii (see~\cite{LyndonSchupp}), which we shall briefly describe in what follows, consider two cases:

Case 1. There exists a letter $t\in S$ whose exponent sum in $r$ is $0$.
To fix a notation, let $S=\{ t, b, c, d, \dots \}$, and, by means of a cyclic permutation or $r$, one may assume that the latter word begins with $b$ or $b^{-1}$. 

Let $b_i$ denote $t^ibt^{-i}$ for $i\in\mathbb Z$, $c_i$ denote $t^ict^{-i}$ for $i\in\mathbb Z$, and so on.
Rewrite $r$ scanning it from left to right and changing any occurrence of $t^ix$ into $t^ixt^{-i}t^i=x_it^i$ (here $x$ represents any letter among $b, c, d, \dots$, or their inverses), collecting the powers of adjacent $t$-letters together, and continuing with the leftmost occurrence of $t$ or its inverse in the modified word. This way we do at least one cancellation of $t$ and its inverse which happen to be next to each other, and the resulting word $s$, which represents $r$ in terms of $t$, $b_i$, $c_i$, \dots, and their inverses, has length not more than $(|r|-2)$. 

Let $m$ and $M$ be the minimal and the maximal subscript of $b_i$ occurring in $s$. Then
\begin{multline*}G\cong< t, b_m, \dots, b_M, c_i, d_i, \dots (i\in\mathbb Z) | \\
s, tb_it^{-1}b_{i+1}^{-1} (i=m, \dots, M-1), tc_it^{-1}c_{i+1}^{-1}, td_it^{-1}d_{i+1}^{-1}, \dots (i\in\mathbb Z)>.\end{multline*}
Consider
$$H=< b_m, \dots, b_M, c_i, d_i, \dots (i\in\mathbb Z) | s >.$$
According to our inductive assumption, $\asdim H\leq \lceil\frac{|s|}2\rceil\leq\lceil\frac{|r|}2\rceil-1$. Now $G\cong*_{F<b_m, \dots, b_M, c_i, d_i, \dots (i\in\mathbb Z)>}H$, and, via Theorem~\ref{thmHNNextension}, $\asdim G\leq 1+\asdim H\leq\lceil\frac{|r|}2\rceil$.

Case 2. For \emph{all} letters in $S$, their exponent sums in $r$ are nonzero. Let $S=\{ u, v, c, d, \dots \}$, and assume that the exponent sums of $u$ and $v$ in $r$ are $\alpha$ and $\beta$ respectively. Define the following homomorphism:
$$\Psi: u\mapsto bt^{-\beta}, v\mapsto t^{\alpha}, c\mapsto c, d\mapsto d, \dots$$
Our group $G$ embeds via $\Psi$ into
$$C=<t, b, c, d, \dots | r( bt^{-\beta}, t^{\alpha}, c, d, \dots) >,$$
and if $p$ is the cyclically reduced $r( bt^{-\beta}, t^{\alpha}, c, d, \dots)$, the exponent sum of $t$ in $p$ is $0$, and $b$ occurs in $p$.

Since $G$ can be thought of as a subgroup of $C$, it would be enough to show that $\asdim C\leq\lceil\frac{|r|}2\rceil$.

Now $C$ is an HNN extension of some group $H$ as in Case 1: assuming that $p$ starts with $b$ or $b^{-1}$, we introduce new variables $b_i=t^ibt^{-i}$, $c_i=t^ict^{-i}$, and so on for $i\in\mathbb Z$. Using these variables, we rewrite $p$ (and therefore $r$ as well) as a word $s$, eliminating all occurrences of $t$ and its inverse and substituting appropriate $x_i$ for any other letter $x$ among $b, c, d,\dots$ and their inverses. If $p$ had at least two occurrences of $t$ or $t^{-1}$, then $|s|\leq|r|-2$, and, using our inductive assumption for $s$, $\asdim H\leq \lceil\frac{|r|-2}2\rceil$ as before. Invoking Theorem~\ref{thmHNNextension}, we obtain the desired inequality $\asdim C\leq\lceil\frac{|r|}2\rceil$.

If, however, $p$ contains $t$ or its inverse only at one place, $p$ expresses $t$ in terms of other generators, so that we can eliminate $t$ from the generating set, and $C$ is indeed a free group on remaining generators with $\asdim C=1\leq\lceil\frac{|r|}2\rceil$.

Now the theorem is proven completely.


\end{document}